\documentclass[12pt,twoside,doublespace]{article}

\usepackage{amssymb}
\usepackage{graphicx}

\def\smskip{\par\vskip 5 pt}
\def\QED{\hfill $\Box$\smskip}
\newtheorem{theorem}{Theorem}
\newtheorem{lemma}{Lemma}
\newtheorem{proposition}{Proposition}

\begin{document}

\begin{center}

\vspace{35pt}

{\Large \bf Gradient Methods with Regularization}

\vspace{5pt}

{\Large \bf  for Constrained Optimization Problems}

\vspace{5pt}

{\Large \bf and Their Complexity Estimates }

\vspace{35pt}

{\sc I.V.~Konnov\footnote{\normalsize E-mail: konn-igor@ya.ru}}

\vspace{35pt}

{\em  Department of System Analysis
and Information Technologies, \\ Kazan Federal University, ul.
Kremlevskaya, 18, Kazan 420008, Russia.}

\end{center}

\vspace{35pt}

\begin{abstract}
We suggest simple implementable modifications of conditional gradient
and gradient projection methods for smooth convex optimization problems in Hilbert
spaces. Usually, the custom methods attain only weak convergence.
We prove strong convergence of the new versions and establish their complexity
estimates, which appear similar to the convergence rate of the
weakly convergent versions.

\end{abstract}


\newpage


\section{Introduction} \label{sc:1}

Let $D$ be a  convex set in a real Hilbert space $H$ and $f : D
\rightarrow \mathbb{R}$  a convex function. Then one can define the
optimization problem of finding the minimal value of the function
$f$ over the feasible set  $D$. For brevity, we write this problem as
\begin{equation} \label{eq:1.1}
 \min \limits _{x \in D} \to f(x),
\end{equation}
its solution set is denoted by $D^{*}(f)$ and the optimal value of
the function by $f^{*}$, i.e.
$$
f^{*} = \inf \limits _{ x \in D} f(x).
$$
For many significant applications this problem appears ill-posed, i.e.
its solution does not depend continuously on the input data. At the same time,
the custom convex optimization methods can in general
provide only weak convergence to a solution, hence, they do not guarantee
sufficient distance approximation of the solution set $D^{*}(f)$, besides,
 even small perturbation of the input data may give large deviations from
the solution. In order to overcome these
drawbacks, various regularization techniques that yield the strong
convergence can be applied; see e.g.
\cite{Vas81}--\cite{EHN96}. The most popular and efficient regularization
 method was suggested by A.N.~Tikhonov; see \cite{Tik63}.

That is, a family of perturbed problems with better properties is solved instead of the
initial one. However, the solution of such a perturbed problem
within a prescribed accuracy may be too difficult even for the
convex optimization problem (\ref{eq:1.1}). At the same time, various
simple and implementable versions of the regularization methods
yield slow convergence due to the special restrictive rules for the choice of step-size
and regularization parameters; see e.g. \cite{BG89,VA93}.

In this paper, we suggest an intermediate variant of the implementable
regularization method. We take the conditional gradient and gradient projection methods as
basic ones. At each iteration of the selected method
it is applied to some perturbed convex
optimization problem. Unlike the known iterative regularization
methods (see \cite{BG89}), we change the perturbed problem only after satisfying some
simple estimate inequality, which allows us to utilize rather mild
rules for the choice of the parameters. Within these rules we prove strong
convergence and establish some complexity estimates for these two-level methods. In
particular, they show that this way of incorporating the
regularization techniques gives almost the same convergence rate as
the custom single-level methods, which provide only
weak convergence.


\section{Properties of regularization
methods} \label{sc:2}

We first recall some definitions. Given a set $X$, a function $f : X
\rightarrow \mathbb{R}$ is said to be

(a) {\em convex}, if for each pair of points $x, y \in X$ and for
all $\alpha \in [0, 1]$, it holds that
$$
f (\alpha x+(1-\alpha )y)\leq \alpha f (x)+(1-\alpha )f (y);
$$

(b) {\em strongly convex} with constant $\varkappa > 0$, if for each
pair of points $x, y \in X$ and for all $\alpha \in [0, 1]$, it
holds that
$$
f (\alpha x+(1-\alpha )y)\leq \alpha f (x)+(1-\alpha )f
(y)-0.5\varkappa\alpha (1-\alpha )\|x-y\|^{2};
$$

(c) {\em  upper (lower) semicontinuous} at a point  $z \in X$, if
for each sequence $\{ x^{k}\} \to z$, $x^{k} \in X$, it holds that
$$
\limsup_{k \to \infty} f (x^{k}) \le f (z) \ (\liminf_{k \to \infty}
f (x^{k}) \ge f (z)).
$$

We will consider problem (\ref{eq:1.1}) under the following basic
assumptions.

 {\bf (A1)} {\em $D$ is a nonempty, convex and closed subset
of a real Hilbert  space $H$, $f : D \rightarrow \mathbb{R}$ is a
lower semicontinuous and convex function.}

The classical Tikhonov regularization method (see \cite{Tik63})
consists in replacing problem (\ref{eq:1.1}) with a sequence of
perturbed problems of the form
\begin{equation} \label{eq:2.1}
 \min \limits _{x \in D} \to  \{f(x)+\varepsilon \varphi(x)\},
\end{equation}
where $\varphi : H \rightarrow \mathbb{R}$ is a lower semicontinuous
and strongly convex function, $\varepsilon>0$ is a regularization
parameter. We recall the basic approximation property; see e.g.
\cite[Chapter II, Section 5, Theorem 1]{Vas81}.


\begin{proposition} \label{pro:2.1} Suppose that all the assumptions in
 {\bf (A1)} are fulfilled,  $D^{*}(f)\neq \varnothing$, and that
 $\varphi : H \rightarrow \mathbb{R}$ is a lower semicontinuous
and strongly convex function. Then:

(i) problem {\em (\ref{eq:2.1})} has a unique solution
$z(\varepsilon)$ for each $\varepsilon>0$;

(ii) if $\{\varepsilon_{k} \} \searrow 0$ as $k \to +\infty$, the
corresponding sequence $\{z(\varepsilon_{k}) \}$ converges strongly
to the point $x^{*}_{n}$ that is the unique solution of the problem
$$
\min \limits _{x \in D^{*}(f)} \to \varphi(x).
$$
\end{proposition}

The main issue of the above regularization method consists in its
suitable implementation since we can not find the point
$z(\varepsilon)$ exactly in the general nonlinear case. Clearly,
instead of $z(\varepsilon)$ we can in principle take any point
$\tilde z(\varepsilon) \in D$ such that $\|\tilde
z(\varepsilon)-z(\varepsilon)\| \leq \xi (\varepsilon)$ with $\xi
(\varepsilon) \searrow 0$ as $\varepsilon \searrow 0$. Then
$\{\tilde z(\varepsilon_{k}) \}$ also converges strongly to the point
$x^{*}_{n}$ in case (ii) of Proposition \ref{pro:2.1}. However, it
is not so easy to guarantee even the prescribed distance
approximation to the point $z(\varepsilon)$  in the general
case.

In \cite{BP74}, the so-called iterative regularization method was
proposed; see \cite{BG89} for more details.
The idea of this method consists in simultaneous changes
of the regularization parameters and step-sizes of a chosen basic
approximation method. In particular, if the functions $f : D
\rightarrow \mathbb{R}$ and $\varphi : H \rightarrow \mathbb{R}$ are
smooth, we can take the basic gradient projection method for problem
(\ref{eq:2.1}). Then the corresponding iterative procedure can be
determined as follows:
\begin{equation} \label{eq:2.2}
\displaystyle x^{k+1}=\pi_{D}[x^k-\lambda_k (f'(x^k)+\varepsilon_k
\varphi'(x^k))], \quad \varepsilon_k>0, \ \lambda_k >0,  \ k=0,1,\ldots;
\end{equation}
where
\begin{equation} \label{eq:2.3}
\begin{array}{lc}
\displaystyle \lim_{k\rightarrow \infty } \varepsilon_k=0, \
  \lim_{k\rightarrow \infty } (\lambda_k/\varepsilon_k ) =0, &\\
 \displaystyle
  \lim_{k\rightarrow \infty } \frac{\varepsilon_k-\varepsilon_{k+1} }
        {\lambda_k\varepsilon_k^2 } =0, \
  \sum \limits_{k=0}^{\infty } (\varepsilon_k\lambda_k)=\infty; & \\
 \end{array}
\end{equation}
and $x^{0} \in D$. Here and below, $\pi _{X} (x)$ denotes the
projection of $x$ onto $X$.


\begin{proposition} \label{pro:2.2} \cite[Theorem 3.1]{BG89}
Suppose that all the assumptions in
 {\bf (A1)} are fulfilled,  $D^{*}(f)\neq \varnothing$, the function $f : D
\rightarrow \mathbb{R}$ is smooth, the function $\varphi : H \rightarrow
\mathbb{R}$ is smooth and strongly convex,  there exists a
constant $M$ such that
$$
\|f'(x)\|\leq M(1+\|x\|) \ \mbox{and} \ \|\varphi'(x)\|\leq
M(1+\|x\|) \quad \forall x\in D.
$$
Then any sequence $\{x^{k}\}$ generated in conformity with rules
(\ref{eq:2.2}) -- (\ref{eq:2.3}) converges strongly to the point
$x^*_{n}$.
\end{proposition}

Of course, the implementation of method (\ref{eq:2.2}) --
(\ref{eq:2.3}) is relatively simple.  Observe that the conditions in
(\ref{eq:2.3}) are fulfilled if we set
$$
\lambda _{k}=(k+1)^{-0.5}, \varepsilon _{k}=(k+1)^{-\tau }, \tau \in
(0, 0.5).
$$
This means that the convergence of the iterative regularization
method may be rather slow in comparison with that of the basic
method. In fact, let us consider the custom gradient projection
method:
\begin{equation} \label{eq:2.4}
\displaystyle x^{k+1}=\pi_{D}[x^k-\lambda_k f'(x^k)],  \ \lambda_k
>0, \ k=0,1,\ldots,
\end{equation}
and $x^{0} \in D$. For brevity, set $\Delta(x)=f(x)-f^{*}$.

{\bf (A2)} {\em The function $f : D \rightarrow \mathbb{R}$ is
smooth and its gradient  satisfies the Lipschitz condition with
constant $L$.}


\begin{proposition} \label{pro:2.3} (\cite[Theorem 5.1]{LP66}
and \cite[Chapter III, Theorem 2.6]{DR68}) Suppose that {\bf (A1)} and {\bf (A2)}
are fulfilled, a sequence $\{x^{k}\}$ is generated in conformity
with rule (\ref{eq:2.4}) where
\begin{equation} \label{eq:2.5}
\lambda_k \in [\lambda',\lambda''], \ \lambda'
>0, \lambda''<2/L.
\end{equation}
Then these exists some constant $C<+\infty$ such that
\begin{equation} \label{eq:2.6}
\Delta(x^k)\leq C/k \quad \mbox{for} \ k=0,1,\ldots
\end{equation}
\end{proposition}
It is well known that method (\ref{eq:2.4}) -- (\ref{eq:2.5}),
unlike (\ref{eq:2.2}) -- (\ref{eq:2.3}), provides only weak
convergence. At the same time, comparing the step-size rules
(\ref{eq:2.3}) and (\ref{eq:2.5}) we can conclude that it seems
rather difficult to obtain the estimate similar to (\ref{eq:2.6})
for the iterative regularization method (\ref{eq:2.2}) --
(\ref{eq:2.3}). The same convergence properties were established for
the gradient projection method with some other known step-size rules
such as the exact one-dimensional minimization and Armijo rules.


\section{Two-level gradient projection method with regularization}\label{sc:3}

We now describe some other way to create an implementable
regularization method, which is based on the gradient projection
method. The method is applied to problem (\ref{eq:1.1}) under the
assumptions {\bf (A1)} and {\bf (A2)}.  At each iteration, the
gradient projection method is applied to some perturbed problem of
form (\ref{eq:2.1}), however, the perturbed problem is changed only
after satisfying some simple estimate inequality, unlike the above
regularization methods. For the simplicity of exposition, we take
the standard perturbation function $\varphi(x)=0.5\|x\|^{2}$, then
we rewrite the perturbed problem
\begin{equation} \label{eq:3.1}
\min \limits _{x \in D} \to \varphi_{\varepsilon}(x)=
\{f(x)+0.5\varepsilon \|x\|^{2}\},
\end{equation}
and set
$$
\varphi^{*}_{\varepsilon} = \inf \limits _{ x \in D}
\varphi_{\varepsilon}(x).
$$
Observe that problem  (\ref{eq:3.1}) has the unique solution
$z(\varepsilon)$ for each $\varepsilon>0$ under the
assumptions {\bf (A1)} and {\bf (A2)} due to Proposition
\ref{pro:2.1} (i), hence $\varphi^{*}_{\varepsilon}
=\varphi_{\varepsilon}(z(\varepsilon))$. Denote by $\mathbb{Z}_{+}$
the set of non-negative integers.

\medskip
\noindent {\bf Method (GPRM).}

 {\em Step 0:} Choose a point
$w^{0}\in D$, numbers $\beta \in (0,1)$, $\theta \in (0,1)$,
sequences $\{\delta _{l}\} \searrow 0$ and $\{\varepsilon_{l}\}
\searrow 0$. Set $l=1$.

{\em Step 1:} Set $x^{0}=w^{l-1}$, $k=0$.

{\em Step 2:}  Take
$y^{k}=\pi_{D}[x^k-\varphi_{\varepsilon_{l}}'(x^k)]$. If
\begin{equation} \label{eq:3.2}
 \|x^{k}-y^{k}\| \leq \delta _{l},
\end{equation}
 set $w^{l}={\rm argmin}\{\varphi_{\varepsilon_{l}}(x^{k}),\varphi_{\varepsilon_{l}}(y^{k})\}$,
 $l=l+1$ and go to Step 1.
{\em (Change the perturbation)}

{\em Step 3:}  Set $d^{k}=y^{k}-x^{k}$, determine $m$ as the
smallest number in $\mathbb{Z}_{+}$ such that
\begin{equation} \label{eq:3.3}
 \varphi_{\varepsilon_{l}}(x^{k}+\theta ^{m} d^{k})
 \leq \varphi_{\varepsilon_{l}} (x^{k})-\beta \theta ^{m}\|d^{k}\|^{2},
\end{equation}
set $\lambda_{k}=\theta ^{m}$, $x^{k+1}=x^{k}+\lambda_{k}d^{k}$,
$k=k+1$, and go to Step 2.
\medskip

We see that the upper level changes the current perturbed problem
which is associated to the index $l$, whereas the lower level with iterations in $k$ is
nothing but the custom gradient projection method with the Armijo
step-size rule applied to the fixed perturbed problem (\ref{eq:3.1})
with $\varepsilon=\varepsilon_{l}$. Clearly, condition
(\ref{eq:3.2}) is very simple and suitable for the verification.

We now give some useful properties of the gradient projection
method.


\begin{lemma} \label{lm:3.1} Suppose that {\bf (A1)} and {\bf (A2)}
are fulfilled. Fix any  $l$. Then we have
\begin{equation} \label{eq:3.4}
  \langle  \varphi_{\varepsilon_{l}}'(x^k)+y^{k}-x^{k}, x-y^{k} \rangle    \geq 0
      \quad \forall x\in D;
\end{equation}
for any $k=0,1,\ldots $; besides, $\lambda_{k}
\geq  \gamma>0$ for any $k=0,1,\ldots $
\end{lemma}
{\bf Proof.}
Relation (\ref{eq:3.4}) follows directly for the projection
properties. Next, under the assumptions made the gradient of the
function $\varphi_{\varepsilon_{l}}$  satisfies the Lipschitz
condition with constant $L'=L+\varepsilon_{0}$. Hence, for any pair
of points $x,y$ we now have
$$
\varphi_{\varepsilon_{l}} (y) \leq
\varphi_{\varepsilon_{l}}(x)+\langle
\varphi_{\varepsilon_{l}}'(x),y-x \rangle +0.5L' \|y-x \|^{2};
$$
see \cite[Chapter III, Lemma 1.2]{DR68}. Then  (\ref{eq:3.4}) gives
\begin{eqnarray*}
\displaystyle && \varphi_{\varepsilon_{l}}(x^{k}+\lambda d^{k})
 - \varphi_{\varepsilon_{l}}(x^{k})
   \leq \lambda \langle \varphi_{\varepsilon_{l}}'(x^{k}),d^{k} \rangle
          +0.5L' \lambda^{2}\| d^{k} \|^{2} \\
 && \leq -\lambda (1-0.5L' \lambda) \| d^{k} \|^{2}  \leq -\beta \lambda \|d^{k} \|^{2},
\end{eqnarray*}
if $ \lambda \leq \bar \lambda=2(1-\beta) /L'$. It follows from
(\ref{eq:3.3}) that $\lambda_{k} \geq \gamma= \min\{1, \theta \bar
\lambda\}>0$.
\QED

We  show that the sequence of perturbed problems is infinite.


\begin{lemma} \label{lm:3.2} Suppose that {\bf (A1)} and {\bf (A2)}
are fulfilled. Then the number of iterations in $k$ for each number
$l$ is finite.
\end{lemma}
{\bf Proof.}
It follows from (\ref{eq:3.3}) and Lemma \ref{lm:3.1} that
$\varphi_{\varepsilon_{l}}(x^{k+1})\leq
\varphi_{\varepsilon_{l}}(x^{k})-\beta \gamma \| d^{k} \|^{2}$, but
$\varphi^{*}_{\varepsilon} >- \infty$, hence $\lim
\limits_{k\rightarrow \infty }d^{k}=\mathbf{0}$, and the result
follows.
\QED

The next property enables us to evaluate the approximation error.


\begin{lemma} \label{lm:3.3} Suppose that {\bf (A1)} and {\bf (A2)}
are fulfilled. Fix any  $l$. Then
\begin{equation} \label{eq:3.5}
 0.5\varepsilon_{l} \| y^{k}- z(\varepsilon_{l})\|^{2} \leq
 \varphi_{\varepsilon_{l}}(y^{k})-\varphi^{*}_{\varepsilon_{l}} \leq
 (L'+1)\| y^{k}-x^{k}\|\| y^{k}- z(\varepsilon_{l})\|
\end{equation}
for any $k=0,1,\ldots $
\end{lemma}
{\bf Proof.} Since $\varphi_{\varepsilon_{l}}$ is strongly
convex with modulus $\varepsilon_{l}$, we have
$$
0.5\varepsilon_{l} \| y^{k}- z(\varepsilon_{l})\|^{2} \leq
 \varphi_{\varepsilon_{l}}(y^{k})-\varphi^{*}_{\varepsilon_{l}} \leq
 \langle \varphi_{\varepsilon_{l}}'(y^{k}),y^{k}-z(\varepsilon_{l})
 \rangle ;
$$
see e.g. \cite[Chapter I, Section 2]{Vas81}. Next, (\ref{eq:3.4})
gives
\begin{eqnarray*}
\displaystyle && \langle
\varphi_{\varepsilon_{l}}'(y^{k}),y^{k}-z(\varepsilon_{l})
 \rangle  \leq \langle \varphi_{\varepsilon_{l}}'(y^{k})-\varphi_{\varepsilon_{l}}'(x^{k})-(y^{k}-x^{k}),
 y^{k}-z(\varepsilon_{l}) \rangle \\
&&
+\langle\varphi_{\varepsilon_{l}}'(x^{k})+(y^{k}-x^{k}),y^{k}-z(\varepsilon_{l})
 \rangle \\
 && \leq \langle \varphi_{\varepsilon_{l}}'(y^{k})-\varphi_{\varepsilon_{l}}'(x^{k})-(y^{k}-x^{k}),
 y^{k}-z(\varepsilon_{l}) \rangle \\
 && \leq (\| \varphi_{\varepsilon_{l}}'(y^{k})-\varphi_{\varepsilon_{l}}'(x^{k}) \|
  +\| y^{k}-x^{k}\|)\| y^{k}- z(\varepsilon_{l})\|  \\
 &&  \leq (L'+1)\| y^{k}-x^{k}\|\| y^{k}- z(\varepsilon_{l})\|.
\end{eqnarray*}
It follows that (\ref{eq:3.5}) holds true.
\QED

We are ready to establish the basic convergence property for (GPRM).


\begin{theorem} \label{thm:3.1} Suppose that {\bf (A1)} and {\bf (A2)}
are fulfilled and  $D^{*}(f)\neq \varnothing$, we apply (GPRM) with
\begin{equation} \label{eq:3.6}
\lim \limits_{k\rightarrow \infty }(\delta _{l}/\varepsilon_{l})=0.
\end{equation}
Then:

(i) the number of iterations in $k$ for each number $l$ is finite;

(ii) the  sequence $\{w^{l} \}$ converges strongly to the point
$x^{*}_{n}$.
\end{theorem}
{\bf Proof.}
Assertion (i) has been obtained in Lemma \ref{lm:3.2}. Fix any $l$ and denote by $k(l)$
the maximal value of the index $k$ for this $l$, i.e. $\| y^{k(l)}- x^{k(l)})\| \leq \delta
 _{l}$. Then (\ref{eq:3.5}) gives
$$
 \| y^{k(l)}- z(\varepsilon_{l})\| \leq 2(L'+1)\delta _{l}/\varepsilon_{l},
$$
but
$$
 \|w^{l}- y^{k(l)}+ y^{k(l)}- z(\varepsilon_{l})\| \leq \| w^{l}- y^{k(l)}\| + \|y^{k(l)}- z(\varepsilon_{l})\|
 \leq \delta _{l}+\|y^{k(l)}- z(\varepsilon_{l})\|,
$$
hence
$$
 \| w^{l}- z(\varepsilon_{l})\| \leq \left(2(L'+1)/\varepsilon_{l}+1\right)\delta
 _{l}.
$$
Therefore, by (\ref{eq:3.6}),
$$
\lim \limits_{l\rightarrow \infty }\| w^{l}- z(\varepsilon_{l})\|=0.
$$
Due to Proposition \ref{pro:2.1} (ii), $\{z(\varepsilon_{l}) \}$
converges strongly to $x^{*}_{n}$. Therefore, assertion (ii) is also
true.
\QED

We observe that inserting the control sequence $\{\delta _{l}\}$
does not require additional computational expenses per iteration,
but implies the strong convergence, whereas the usual gradient
projection method provides only weak
convergence as indicated above. Besides, rule (\ref{eq:3.6})
is clearly less restrictive
than (\ref{eq:2.3}) and maintains significant freedom
for the choice of parameters.


\section{Complexity estimate}\label{sc:4}

It was observed in Section \ref{sc:2} that the usual gradient
projection method has the convergence rate $\Delta(x^k)\leq C/k$ under the
assumptions {\bf (A1)} and {\bf (A2)}; see
Proposition \ref{pro:2.3} and the remarks below. This means that
the total number of iterations $N(\alpha)$ that is necessary for attaining some prescribed
accuracy  $\alpha>0$ is estimated as follows:
\begin{equation} \label{eq:4.1}
N(\alpha)\leq C/\alpha.
\end{equation}
We intend to obtain a similar estimate for (GPRM). Namely,
we define the complexity of (GPRM), denoted by $N(\alpha)$,  as the
total number of iterations in $k$
that is necessary for attaining any
accuracy  $\alpha>0$. In order to establish an upper bound for
$N(\alpha)$ we need certain auxiliary properties.
We recall that $z(\varepsilon)$ denotes the solution of the
perturbed problem (\ref{eq:3.1}) for $\varepsilon>0$, which is defined uniquely
under {\bf (A1)}. Hence $z(0)$ denotes any solution of problem (\ref{eq:1.1}).


\begin{lemma} \label{lm:4.1} Suppose that {\bf (A1)} holds.
Then for any numbers $\mu$ and $\eta$ such that $0 \leq \mu <\eta$
we have
\begin{eqnarray}
  && f(z(\eta))-f(z(\mu)) \leq 0.5\eta(\| z(\mu)\|^{2}-\| z(\eta)\|^{2}), \label{eq:4.2}\\
  && \varphi^{*}_{\eta}-\varphi^{*}_{\mu}\leq 0.5(\eta-\mu)\| z(\mu)\|^{2}, \label{eq:4.3}\\
  && \|z(\eta)\| \leq \| z(\mu)\|. \label{eq:4.4}
\end{eqnarray}
\end{lemma}
{\bf Proof.}
By definition,
\begin{eqnarray*}
  \varphi^{*}_{\eta} &=& f(z(\eta))+0.5\eta \| z(\eta)\|^{2}
        \leq f(z(\mu))+0.5\eta \| z(\mu)\|^{2}\\
 &=& f(z(\mu))+0.5\mu \| z(\mu)\|^{2}+0.5(\eta-\mu)\| z(\mu)\|^{2} \\
 &=&\varphi^{*}_{\mu}+0.5(\eta-\mu)\| z(\mu)\|^{2} \\
 & \leq & f(z(\eta))+0.5\mu \| z(\eta)\|^{2} +0.5(\eta-\mu)\| z(\mu)\|^{2}.
\end{eqnarray*}
These relations give (\ref{eq:4.2}) and (\ref{eq:4.3}), besides, we also have
$$
f(z(\eta))+0.5\eta \| z(\eta)\|^{2}
        \leq f(z(\eta))+0.5\mu \| z(\eta)\|^{2} +0.5(\eta-\mu)\| z(\mu)\|^{2},
$$
which gives (\ref{eq:4.4}).
\QED

Denote by $N_{(l)}$ the total number of iterations in $k$ for any fixed
$l$ in (GPRM) and by $l(\alpha )$ the maximal
number $l$ of the upper iteration such that
$\alpha \leq \Delta(w^l)$ for any given
$\alpha>0$. Then we can evaluate
the complexity of (GPRM) as follows:
\begin{equation} \label{eq:4.5}
 N (\alpha ) \leq  \sum ^{l(\alpha  )} _{l=1}
     N_{(l)}.
\end{equation}
Using this inequality, we now obtain the basic estimate.


\begin{theorem} \label{thm:4.1}
Suppose that {\bf (A1)} and {\bf (A2)}
are fulfilled and  $D^{*}(f)\neq \varnothing$, we apply (GPRM) with
\begin{equation} \label{eq:4.6}
 \varepsilon_{l} =\nu ^{l}\varepsilon_{0}, \ \delta _{l} = \varepsilon_{l} ^{1+\sigma},  \ l=0,1,\ldots;
  \quad \nu \in (0,1), \ \sigma \in (0,1], \ \varepsilon_{0}>0.
\end{equation}
 Then (GPRM) has the complexity
estimate
$$
N (\alpha ) \leq C_{2} ((C_{1}/\alpha)^{1+2\sigma}-1)/(\nu(1-\nu^{1+2\sigma})),
$$
where $C_{1}=2 (L'+1)^{2} \varepsilon_{0}^{1+2\sigma}+ 0.5 \varepsilon_{0} \| x^{*}_{n}\|^{2}$ and
$C_{2}=C_{1}/(\beta \gamma \varepsilon^{2(1+\sigma)}_{0})$.
\end{theorem}
{\bf Proof.}
First we note that (\ref{eq:4.6}) implies (\ref{eq:3.6}), hence all the
assertions of Theorem \ref{thm:3.1} remain true.  Fix any $l$.
Then, due to (\ref{eq:3.3})
and Lemma \ref{lm:3.1}, we have
$$
 \varphi_{\varepsilon_{l}}(x^{k+1})
 \leq \varphi_{\varepsilon_{l}}(x^{k}) -\beta \gamma \delta _{l}^{2},
$$
therefore,
\begin{equation} \label{eq:4.7}
 N_{(l)} \leq \left(\varphi_{\varepsilon_{l}}(w^{l-1})-
\varphi^{*}_{\varepsilon_{l}}\right)/(\beta \gamma \delta _{l}^{2}).
\end{equation}
However,
\begin{eqnarray*}
 \varphi_{\varepsilon_{l}}(w^{l-1})- \varphi^{*}_{\varepsilon_{l}}
   &=& f(w^{l-1})+0.5\varepsilon_{l-1} \|w^{l-1}\|^{2}+0.5(\varepsilon_{l}-\varepsilon_{l-1}) \|w^{l-1}\|^{2}- \varphi^{*}_{\varepsilon_{l}} \\
   & \leq &     \varphi_{\varepsilon_{l-1}}(w^{l-1})- \varphi^{*}_{\varepsilon_{l-1}}
   + \varphi^{*}_{\varepsilon_{l-1}}- \varphi^{*}_{\varepsilon_{l}}.
\end{eqnarray*}
From (\ref{eq:3.5}) we have
$$
0.5\varepsilon_{l} \| y^{k}- z(\varepsilon_{l})\| \leq (L'+1)\| y^{k}-x^{k}\|,
$$
hence
\begin{eqnarray*}
 \varphi_{\varepsilon_{l}}( y^{k})- \varphi^{*}_{\varepsilon_{l}}
   &\leq &   (L'+1)\| y^{k}-x^{k}\|\| y^{k}- z(\varepsilon_{l})\| \leq 2 ((L'+1)\| y^{k}-x^{k}\|)^{2}/\varepsilon_{l}.
\end{eqnarray*}
It follows that
\begin{equation} \label{eq:4.8}
 \varphi_{\varepsilon_{l}}(w^{l})- \varphi^{*}_{\varepsilon_{l}} \leq   2 ((L'+1) \delta _{l})^{2}/\varepsilon_{l},
\end{equation}
whereas (\ref{eq:4.3}) and (\ref{eq:4.4}) give
$$
 \varphi^{*}_{\varepsilon_{l-1}}- \varphi^{*}_{\varepsilon_{l}}
  \leq   0.5 (\varepsilon_{l-1}-\varepsilon_{l}) \| z(\varepsilon_{l})\|^{2}
   \leq   0.5 (\varepsilon_{l-1}-\varepsilon_{l}) \|x^{*}_{n}\|^{2}.
$$
Therefore,
\begin{eqnarray*}
 \varphi_{\varepsilon_{l}}(w^{l-1})- \varphi^{*}_{\varepsilon_{l}}
   & \leq & 2 ((L'+1) \delta _{l-1})^{2}/\varepsilon_{l-1}+ 0.5 (\varepsilon_{l-1}-\varepsilon_{l}) \|x^{*}_{n}\|^{2} \\
&=& 2 ((L'+1) \delta _{0})^{2}\nu^{2(1+\sigma)(l-1)}/(\varepsilon_{0}\nu^{l-1})+ 0.5 (1-\nu)\varepsilon_{0}\nu^{l-1} \|x^{*}_{n}\|^{2} \\
& \leq & C_{1}\nu^{l-1},
\end{eqnarray*}
where
$$
C_{1}=2 (L'+1)^{2} \varepsilon_{0}^{1+2\sigma}+ 0.5 \varepsilon_{0}\|x^{*}_{n}\|^{2}.
$$
Using these relations in (\ref{eq:4.7}) we have
\begin{equation} \label{eq:4.9}
 N_{(l)} \leq (C_{1}\nu^{l})/(\beta \gamma \delta^{2} _{0} \nu^{2(1+\sigma)l+1})=(C_{2}/\nu)\nu^{-(1+2\sigma)l},
\end{equation}
where
$$
C_{2}=C_{1}/(\beta \gamma \varepsilon^{2(1+\sigma)}_{0}).
$$
In view of (\ref{eq:4.5}) and (\ref{eq:4.9}) we obtain
\begin{eqnarray}
 N (\alpha ) & \leq &  (C_{2}/\nu)\sum ^{l(\alpha  )} _{l=1} \nu^{-(1+2\sigma)l}
 =C_{2} \nu^{-(1+2\sigma)-1}\left(\nu^{-(1+2\sigma)l(\alpha  )} -1 \right)/(\nu^{-(1+2\sigma)}-1) \nonumber \\
&=&
C_{2} \left(\nu^{-(1+2\sigma)l(\alpha  )} -1 \right)/(\nu(1-\nu^{1+2\sigma})). \label{eq:4.10}
\end{eqnarray}

We now proceed to evaluate $\nu^{-l(\alpha  )}$. By definition,
\begin{eqnarray*}
 \varphi_{\varepsilon_{l}}(w^{l})- \varphi^{*}_{\varepsilon_{l}}
   &=& f(w^{l})+0.5\varepsilon_{l} \| w^{l}\|^{2} -[ f(z(\varepsilon_{l}))+0.5\varepsilon_{l} \|z(\varepsilon_{l})\|^{2}] \\
 &=& [ f(w^{l})-f^{*}]-[ f(z(\varepsilon_{l}))-f^{*}]
 -0.5\varepsilon_{l}[ \|z(\varepsilon_{l})\|^{2}-\| w^{l}\|^{2}],
\end{eqnarray*}
hence,
$$
\Delta(w^l)=f(w^{l})-f^{*}=\varphi_{\varepsilon_{l}}(w^{l})- \varphi^{*}_{\varepsilon_{l}}
+[ f(z(\varepsilon_{l}))-f^{*}]
 +0.5\varepsilon_{l}[ \|z(\varepsilon_{l})\|^{2}-\| w^{l}\|^{2}].
$$
From (\ref{eq:4.8}) we have
$$
\varphi_{\varepsilon_{l}}(w^{l})- \varphi^{*}_{\varepsilon_{l}} \leq 2 ((L'+1) \delta _{l})^{2}/\varepsilon_{l},
$$
whereas applying (\ref{eq:4.2}) with $\mu=0$ and $\eta=\varepsilon_{l}$ gives
$$
 f(z(\varepsilon_{l}))-f^{*}
  \leq   0.5 \varepsilon_{l} (\| x^{*}_{n}\|^{2}-\| z(\varepsilon_{l})\|^{2} ).
$$
Therefore,
$$
\Delta(w^l) \leq  2 ((L'+1) \delta _{l})^{2}/\varepsilon_{l}+ 0.5 \varepsilon_{l} \| x^{*}_{n}\|^{2}.
$$
In view of (\ref{eq:4.6})  we have
$$
\Delta(w^l) \leq  2 (L'+1)^{2} \varepsilon_{0}^{1+2\sigma} \nu^{(1+2\sigma)l}+ 0.5 \varepsilon_{0} \| x^{*}_{n}\|^{2}\nu^{l}
\leq C_{1}\nu^{l}.
$$
It follows that $\nu^{-l(\alpha)} \leq C_{1}/\alpha$. Applying this inequality in
(\ref{eq:4.10})  we obtain
\begin{eqnarray*}
N (\alpha ) &\leq &   C_{2} \left(\nu^{-(1+2\sigma)l(\alpha  )} -1 \right)/(\nu(1-\nu^{1+2\sigma})) \\
 &\leq &   C_{2} \left((C_{1}/\alpha)^{1+2\sigma} -1 \right)/(\nu(1-\nu^{1+2\sigma})),
\end{eqnarray*}
and the result follows.
\QED

From Theorem \ref{thm:4.1} we conclude that the complexity
estimate of (GPRM) tends to (\ref{eq:4.1}) when $\sigma \to 0$.
However, we can choose $\sigma$ arbitrarily in $(0,1]$.
Therefore, taking $\sigma$ small enough, we can obtain any approximation of
the convergence rate of the usual gradient
projection method under the same assumptions. At the same time,
(GPRM), unlike the gradient
projection method, attains the strong convergence.


\section{Two-level conditional gradient method with regularization}\label{sc:5}

We now describe a similar modification of the conditional gradient
method under the following basic
assumptions for problem (\ref{eq:1.1}).

{\bf (A3)} {\em $D$ is a nonempty, convex, closed, and bounded subset
of a real Hilbert  space $H$, $f : D \rightarrow \mathbb{R}$ is
a smooth convex function and its gradient satisfies the Lipschitz condition with
constant $L$.}

The boundedness of $D$ guarantees the method is well-defined. Besides, now
problem (\ref{eq:1.1}) has a solution, i.e. $D^{*}(f)\neq \varnothing$.
We recall that the conditional gradient method
was first suggested in \cite{FW56} for the case when the goal function is quadratic and
the feasible set is polyhedral and further was developed by many authors;
see e.g. \cite{LP66,DR68,PD78,Dun80,Kon17}. The main
idea of this method consists in linearization of the goal function, so that
solution of the linearized problem over the initial feasible set serves for finding the
descent direction.

Following \cite{LP66,DR68}, we describe one of the various versions of
the custom conditional gradient method.

\medskip
\noindent {\bf Method (CGM).}

 {\em Step 0:} Choose a point
$x^{0}\in D$,  set $k=0$.

{\em Step 1:}  Find a point $y^{k}\in D$ as a solution of the problem
$$
 \min \limits _{y \in D} \to \langle  f'(x^k), y\rangle ,
$$
set $d^{k}=y^{k}-x^{k}$.

{\em Step 2:}  If $d^{k}=\mathbf{0}$, stop. Otherwise choose a number
$\theta_{k} >0$, set $\beta_{k} =-\langle  f'(x^k), d^k\rangle/\|d^{k}\|^{2}$,
$\lambda_{k}=\min\{1,\theta_{k}\beta_{k}\}$, $x^{k+1}=x^{k}+\lambda_{k}d^{k}$,
$k=k+1$, and go to Step 1.
\medskip

Clearly, termination of the method yields a solution. For this reason,
we will consider only the non-trivial case where the sequence $\{x^{k}\}$ is infinite.


\begin{proposition} \label{pro:5.1} (\cite[Theorem 6.1]{LP66}
and \cite[Chapter III, Theorem 1.7]{DR68}) Suppose that {\bf (A3)}
is fulfilled, a sequence $\{x^{k}\}$ is generated by (CGM) where
$$
\theta_{k} \in [\theta',\theta''], \ \theta'
>0, \theta''<2/L.
$$
Then these exists some constant $C<+\infty$ such that
\begin{equation} \label{eq:5.1}
\Delta(x^k)\leq C/k \quad \mbox{for} \ k=0,1,\ldots
\end{equation}
\end{proposition}
That is, estimate (\ref{eq:5.1}) is the same as (\ref{eq:2.6}), but
it can not be enhanced even if the function $f$ is strongly convex. Besides,
(CGM) also provides only weak convergence. The same convergence properties
were established for
the conditional gradient method with the other known step-size rules
such as the exact one-dimensional minimization and Armijo rules; see
\cite{DR68,PD78,Dun80}.

Some versions of the iterative regularization method
based on the conditional gradient iterations were described in
\cite[Chapter II, Section 11]{Vas81} and
\cite[Chapter IV, Section 1]{BG89}. They provides strong convergence
but utilize the restrictive control rules for
the regularization parameters and step-sizes, which are similar to (\ref{eq:2.3}).
In particular, the version from \cite{BG89} utilizes
the exact one-dimensional minimization for the choice of the
step-size and take the rule
$$
\varepsilon _{k}=(k+1)^{-\tau }, \tau \in
(0, 0.5),
$$
for the regularization parameter. This means that the convergence of the iterative regularization
version may be rather slow in comparison with that of the basic conditional gradient
method.

We now describe some other implementable conditional gradient
method with regularization, which follows the approach given in Section
\ref{sc:3}. That is, the custom
conditional gradient method is applied to some perturbed problem of
form (\ref{eq:2.1}), however, the perturbed problem is changed only
after satisfying some simple estimate inequality. We also take
the standard perturbation function $\varphi(x)=0.5\|x\|^{2}$, hence
we take the perturbed problem (\ref{eq:3.1}), which has the unique solution
$z(\varepsilon)$ for each $\varepsilon>0$ under the
assumptions in {\bf (A3)}.

\medskip
\noindent {\bf Method (CGRM).}

 {\em Step 0:} Choose a point
$w^{0}\in D$, numbers $\beta \in (0,1)$, $\theta \in (0,1)$,
sequences $\{\delta _{l}\} \searrow 0$ and $\{\varepsilon_{l}\}
\searrow 0$. Set $l=1$.

{\em Step 1:} Set $x^{0}=w^{l-1}$, $k=0$.

{\em Step 2:}  Find a point $y^{k}\in D$ as a solution of the problem
$$
 \min \limits _{y \in D} \to \langle  \varphi_{\varepsilon_{l}}'(x^k), y\rangle ,
$$
set $d^{k}=y^{k}-x^{k}$, $\mu_{k,l}=-\langle  \varphi_{\varepsilon_{l}}'(x^k),d^{k}\rangle$. If
\begin{equation} \label{eq:5.2}
 \mu_{k,l} \leq \delta _{l},
\end{equation}
 set $w^{l}=x^{k}$, $l=l+1$ and go to Step 1.
{\em (Change the perturbation)}

{\em Step 3:}  Determine $m$ as the
smallest number in $\mathbb{Z}_{+}$ such that
\begin{equation} \label{eq:5.3}
 \theta ^{m}\mu_{k,l} \leq 1, \ \varphi_{\varepsilon_{l}}(x^{k}+\theta ^{m}\mu_{k,l} d^{k})
 \leq \varphi_{\varepsilon_{l}} (x^{k})-\beta \theta ^{m}\mu_{k,l}^{2},
\end{equation}
set $\lambda_{k}=\theta ^{m}$, $x^{k+1}=x^{k}+\lambda_{k}\mu_{k,l}d^{k}$,
$k=k+1$, and go to Step 2.
\medskip

We see again that the upper level changes the current perturbed problem
associated to the index $l$, whereas the lower level with iterations in $k$ is
nothing but the conditional gradient method with the Armijo
step-size rule applied to the fixed perturbed problem. Clearly, condition
(\ref{eq:5.2}) is very simple and suitable for the verification.

We now give a lower bound for the step-size.


\begin{lemma} \label{lm:5.1} Suppose that {\bf (A3)}
is fulfilled. Fix any $l$. Then
$$
\lambda_{k} \geq  \gamma >0,
$$
for any $k=0,1,\ldots $
\end{lemma}
{\bf Proof.}
It was noticed that,  under the assumptions made the gradient of the
function $\varphi_{\varepsilon_{l}}$  satisfies the Lipschitz
condition with constant $L'=L+\varepsilon_{0}$. Hence, for any pair
of points $x,y$ we now have
$$
\varphi_{\varepsilon_{l}} (y) \leq
\varphi_{\varepsilon_{l}}(x)+\langle
\varphi_{\varepsilon_{l}}'(x),y-x \rangle +0.5L' \|y-x \|^{2}.
$$
 Therefore,
$$
 \varphi_{\varepsilon_{l}}(x^{k}+\lambda \mu_{k,l} d^{k})
 - \varphi_{\varepsilon_{l}}(x^{k})
   \leq \lambda \mu_{k,l}\langle \varphi_{\varepsilon_{l}}'(x^{k}),d^{k} \rangle
          +0.5L' \lambda^{2}\mu_{k,l}^{2}\| d^{k} \|^{2}
  \leq  -\beta \lambda \mu_{k,l}^{2},
$$
if
$$
(1-\beta )\lambda \mu_{k,l}^{2} \geq 0.5L' \lambda^{2}\mu_{k,l}^{2}\| d^{k} \|^{2}
$$
or
$ \lambda \leq \lambda'=2(1-\beta) /L'B^{2}$, where $B$ denotes the diameter
of the set $D$. Fix any point $\bar x \in D$. Then
\begin{eqnarray*}
  \mu_{k,l} &\leq & \langle
        \varphi_{\varepsilon_{l}}'(\bar x),x^{k}-y^{k} \rangle
        + \langle \varphi_{\varepsilon_{l}}'(x^{k})-
             \varphi_{\varepsilon_{l}}'(\bar x),x^{k}-y^{k} \rangle \\
             &\leq & (\| f'(\bar x) \| + \varepsilon_{0} \|\bar x \| ) B+ L'B^{2} = L''B,
\end{eqnarray*}
hence setting $\lambda'' = 1/(L''B)$ gives $\mu_{k,l} \lambda'' \leq 1 $.
Set $\gamma=\min\{\theta, \lambda', \lambda''\}>0$.
It follows now from
(\ref{eq:5.3}) that $\lambda_{k} \geq \gamma$.
\QED

We now show that the sequence of perturbed problems is infinite.


\begin{lemma} \label{lm:5.2} Suppose that {\bf (A3)}
is fulfilled. Then the number of iterations in $k$ for each number
$l$ is finite.
\end{lemma}
{\bf Proof.}
It follows from (\ref{eq:5.3}) that
$\varphi_{\varepsilon_{l}}(x^{k+1})\leq
\varphi_{\varepsilon_{l}}(x^{k})-\beta \gamma \mu_{k,l}^{2}$, but
$\varphi^{*}_{\varepsilon} >- \infty$, hence $\lim
\limits_{k\rightarrow \infty }\mu_{k,l}=0$, and the result
follows.
\QED

The next property enables us to evaluate the approximation error.


\begin{lemma} \label{lm:5.3} Suppose that {\bf (A3)}
is fulfilled. Fix any  $l$. Then
\begin{equation} \label{eq:5.4}
 0.5\varepsilon_{l} \| x^{k}- z(\varepsilon_{l})\|^{2} \leq
 \varphi_{\varepsilon_{l}}(x^{k})-\varphi^{*}_{\varepsilon_{l}} \leq
 \mu_{k,l}
\end{equation}
for any $k=0,1,\ldots $
\end{lemma}
{\bf Proof.} Since $\varphi_{\varepsilon_{l}}$ is strongly
convex with modulus $\varepsilon_{l}$, we have
$$
0.5\varepsilon_{l} \| x^{k}- z(\varepsilon_{l})\|^{2} \leq
 \varphi_{\varepsilon_{l}}(x^{k})-\varphi^{*}_{\varepsilon_{l}} \leq
 \langle \varphi_{\varepsilon_{l}}'(x^{k}),x^{k}-z(\varepsilon_{l})
 \rangle ;
$$
see e.g. \cite[Chapter I, Section 2]{Vas81}. By definition, we have
$$
\langle \varphi_{\varepsilon_{l}}'(x^{k}),x^{k}-z(\varepsilon_{l})
 \rangle = \langle \varphi_{\varepsilon_{l}}'(x^{k}),x^{k}-y^{k}
 \rangle + \langle \varphi_{\varepsilon_{l}}'(x^{k}),y^{k}-z(\varepsilon_{l})
 \rangle \leq \mu_{k,l},
$$
which gives (\ref{eq:5.4}).
\QED

We are ready to establish the basic convergence property for (CGRM).


\begin{theorem} \label{thm:5.1} Suppose that {\bf (A3)}
is fulfilled, we apply (CGRM) with (\ref{eq:3.6}).
Then:

(i) the number of iterations in $k$ for each number $l$ is finite;

(ii) the  sequence $\{w^{l} \}$ converges strongly to the point
$x^{*}_{n}$.
\end{theorem}
{\bf Proof.}
Assertion (i) has been obtained in Lemma \ref{lm:5.2}. Fix any $l$ and denote by $k(l)$
the maximal value of the index $k$ for this $l$. Then $\mu_{k(l),l} \leq \delta
 _{l}$ and (\ref{eq:5.4}) gives
$$
 \| w^{l}- z(\varepsilon_{l})\|^{2} \leq 2\delta _{l}/\varepsilon_{l},
$$
hence, by (\ref{eq:3.6}),
$$
\lim \limits_{l\rightarrow \infty }\| w^{l}- z(\varepsilon_{l})\|=0.
$$
Due to Proposition \ref{pro:2.1} (ii), $\{z(\varepsilon_{l}) \}$
converges strongly to $x^{*}_{n}$. Therefore, assertion (ii) is also
true.
\QED

We also notice that  rule (\ref{eq:3.6})
is clearly less restrictive
than (\ref{eq:2.3}) and maintains significant freedom
for the choice of parameters.

Due to Proposition \ref{pro:5.1},
the total number of iterations $N(\alpha)$ of the conditional gradient method
that is necessary for attaining some prescribed
accuracy  $\alpha>0$ is estimated as follows:
\begin{equation} \label{eq:5.5}
N(\alpha)\leq C/\alpha.
\end{equation}
We intend to obtain a similar estimate for (CGRM). As above in Section \ref{sc:4},
we define the complexity of (CGRM), denoted by $N(\alpha)$,  as the
total number of iterations in $k$
that is necessary for attaining any given
accuracy  $\alpha>0$.

Denote by $N_{(l)}$ the total number of iterations in $k$ for any fixed
$l$ in (CGRM) and by $l(\alpha )$ the maximal
number $l$ of the upper iteration such that
$\alpha \leq \Delta(w^l)$ for any given
$\alpha>0$. Then we can evaluate
the complexity of (CGRM) as follows:
\begin{equation} \label{eq:5.6}
 N (\alpha ) \leq  \sum ^{l(\alpha  )} _{l=1}
     N_{(l)};
\end{equation}
cf. (\ref{eq:4.5}).
Using this inequality, we now obtain the basic estimate. Its substantiation
is somewhat different from the proof of  Theorem \ref{thm:4.1}.


\begin{theorem} \label{thm:5.2}
Suppose that {\bf (A3)}
is fulfilled, we apply (CGRM) with (\ref{eq:4.6}).
 Then (CGRM) has the complexity
estimate
$$
N (\alpha ) \leq C_{2} ((C_{1}/\alpha)^{1+2\sigma}-1)/(\nu(1-\nu^{1+2\sigma})),
$$
where $C_{1}=\varepsilon_{0}^{1+2\sigma}+ 0.5 \varepsilon_{0}\|x^{*}_{n}\|^{2}$ and
$C_{2}=C_{1}/(\beta \gamma \varepsilon^{2(1+\sigma)}_{0})$.
\end{theorem}
{\bf Proof.}
First we note that (\ref{eq:4.6}) implies (\ref{eq:3.6}), hence all the
assertions of Theorem \ref{thm:5.1} remain true.  Fix any $l$.
Then, due to (\ref{eq:5.3}) and Lemma \ref{lm:5.1}, we have
$$
 \varphi_{\varepsilon_{l}}(x^{k+1})
 \leq \varphi_{\varepsilon_{l}}(x^{k}) -\beta \gamma \delta _{l}^{2},
$$
therefore,
\begin{equation} \label{eq:5.7}
 N_{(l)} \leq \left(\varphi_{\varepsilon_{l}}(w^{l-1})-
\varphi^{*}_{\varepsilon_{l}}\right)/(\beta \gamma \delta _{l}^{2}).
\end{equation}
However,
$$
 \varphi_{\varepsilon_{l}}(w^{l-1})- \varphi^{*}_{\varepsilon_{l}}
    \leq     \varphi_{\varepsilon_{l-1}}(w^{l-1})- \varphi^{*}_{\varepsilon_{l-1}}
   + \varphi^{*}_{\varepsilon_{l-1}}- \varphi^{*}_{\varepsilon_{l}}.
$$
From (\ref{eq:5.4}) we have
\begin{equation} \label{eq:5.8}
 \varphi_{\varepsilon_{l}}(w^{l})- \varphi^{*}_{\varepsilon_{l}} \leq   \delta _{l},
\end{equation}
whereas (\ref{eq:4.3}) and (\ref{eq:4.4}) give
$$
 \varphi^{*}_{\varepsilon_{l-1}}- \varphi^{*}_{\varepsilon_{l}}
  \leq   0.5 (\varepsilon_{l-1}-\varepsilon_{l}) \| z(\varepsilon_{l})\|^{2}
   \leq   0.5 (\varepsilon_{l-1}-\varepsilon_{l}) \|x^{*}_{n}\|^{2}.
$$
Therefore,
$$
 \varphi_{\varepsilon_{l}}(w^{l-1})- \varphi^{*}_{\varepsilon_{l}}
   \leq   \delta _{l-1}+ 0.5 (\varepsilon_{l-1}-\varepsilon_{l}) \|x^{*}_{n}\|^{2}
 \leq  C_{1}\nu^{l-1},
$$
where
$$
C_{1}=\varepsilon_{0}^{1+2\sigma}+ 0.5 \varepsilon_{0}\|x^{*}_{n}\|^{2}.
$$
Using these relations in (\ref{eq:5.7}) we have
\begin{equation} \label{eq:5.9}
 N_{(l)} \leq (C_{1}\nu^{l})/(\beta \gamma \varepsilon^{2(1+\sigma)}_{0} \nu^{2(1+\sigma)l+1})=(C_{2}/\nu)\nu^{-(1+2\sigma)l},
\end{equation}
where
$$
C_{2}=C_{1}/(\beta \gamma \varepsilon^{2(1+\sigma)}_{0}).
$$
In view of (\ref{eq:5.6}) and (\ref{eq:5.9}) we obtain
\begin{eqnarray}
 N (\alpha ) & \leq &  (C_{2}/\nu)\sum ^{l(\alpha  )} _{l=1} \nu^{-(1+2\sigma)l}
 =C_{2} \nu^{-(1+2\sigma)-1}\left(\nu^{-(1+2\sigma)l(\alpha  )} -1 \right)/(\nu^{-(1+2\sigma)}-1) \nonumber \\
&=&
C_{2} \left(\nu^{-(1+2\sigma)l(\alpha  )} -1 \right)/(\nu(1-\nu^{1+2\sigma})). \label{eq:5.10}
\end{eqnarray}

We now proceed to evaluate $\nu^{-l(\alpha  )}$. By definition,
$$
\Delta(w^l)=f(w^{l})-f^{*}=\varphi_{\varepsilon_{l}}(w^{l})- \varphi^{*}_{\varepsilon_{l}}
+[ f(z(\varepsilon_{l}))-f^{*}]
 +0.5\varepsilon_{l}[ \|z(\varepsilon_{l})\|^{2}-\| w^{l}\|^{2}].
$$
Applying (\ref{eq:4.2}) with $\mu=0$ and $\eta=\varepsilon_{l}$ gives
$$
 f(z(\varepsilon_{l}))-f^{*}
  \leq   0.5 \varepsilon_{l} (\| x^{*}_{n}\|^{2}-\| z(\varepsilon_{l})\|^{2} ).
$$
From (\ref{eq:5.8}) it now follows that
$$
\Delta(w^l) \leq  \delta _{l}+ 0.5 \varepsilon_{l} \| x^{*}_{n}\|^{2}.
$$
In view of (\ref{eq:4.6})  we have
$$
\Delta(w^l) \leq   \varepsilon_{0}^{1+2\sigma} \nu^{(1+2\sigma)l}+ 0.5 \varepsilon_{0} \| x^{*}_{n}\|^{2}\nu^{l}
\leq C_{1}\nu^{l}.
$$
It follows that $\nu^{-l(\alpha)} \leq C_{1}/\alpha$. Applying this inequality in
(\ref{eq:5.10})  we obtain
$$
N (\alpha ) \leq   C_{2} \left((C_{1}/\alpha)^{1+2\sigma} -1 \right)/(\nu(1-\nu^{1+2\sigma})),
$$
and the result follows.
\QED

From Theorem \ref{thm:5.2} we conclude that the complexity
estimate of (CGRM) tends to (\ref{eq:5.5}) when $\sigma \to 0$.
Due to (\ref{eq:4.6}), we can choose $\sigma$ arbitrarily in $(0,1]$.
Therefore, taking $\sigma$ small enough, we can obtain any approximation of
the best convergence rate of the usual  conditional gradient
 method under the same assumptions. At the same time,
(CGRM) attains the strong convergence.


\section{Conclusions}\label{sc:6}

We suggested simple implementable versions of the combined
regularization and gradient methods for
smooth convex optimization problems in Hilbert
spaces. We took the basic conditional gradient
and gradient projection methods and
proved strong convergence of their modified versions under rather mild
rules for the choice of the parameters.
Within these rules we also established complexity estimates for the methods.
They show that this way of incorporating the
regularization techniques gives the convergence rate similar to that of
the custom method, which provides only
weak convergence under the same assumptions.


\section*{Acknowledgement}

This work was supported by the RFBR grant, project No.
13-01-00368-a.


\end{document}